\def\IR{{\Bbb R}}

\def\IK{{\Bbb K}}

\documentclass[12pt]{article} 

\usepackage[psamsfonts]{amssymb} 
\usepackage{amsfonts,amsmath} 
\usepackage{tikz-cd}

\usepackage{epsfig,multicol}

\newtheorem{theorem}{Theorem}
\newtheorem{corollary}{Corollary}

\newtheorem{definition}{Definition}[section] 

\title{The Generic Failure of  Lower-semicontinuity for the Linear Distortion Functional}
\author{Sayed Mohsen Hashemi and Gaven J. Martin  \thanks{
Work of both authors partially supported by the New Zealand Marsden Fund.  Parts of this work appear in the PhD thesis of the first author.\newline
{\bf MSC Subject:}  30C60. 
}}

\date{{\em To the memory of Peter Duren}}
\begin{document}
\maketitle
\begin{abstract} {\footnotesize We consider the convexity properties of distortion functionals,  particularly the linear distortion,  defined for homeomorphisms of domains in Euclidean $n$-spaces,  $n\geq 3$. The inner and outer distortion functionals are lower semi-continuous in all dimensions and so for the curve modulus or analytic definitions of quasiconformality it ifollows that if $ \{ f_{n} \}_{n=1}^{\infty} $ is a sequence of $K$-quasiconformal mappings (here $K$ depends on the particular distortion functional but is the same for every element of the sequence) which converges locally uniformly to a mapping $f$, then this limit function is also $K$-quasiconformal.Despite a widespread belief that this was also true for the geometric definition of quasiconformality (defined through the linear distortion $H({f_{n}})$), T. Iwaniec gave a specific and surprising example to show that the linear distortion functional is not always lower semicontinuous on uniformly converging sequences of quasiconformal mappings.  Here we show that this failure of lower semicontinuity is common, perhaps generic in the sense that  under mild restrictions on a quasiconformal $f$, there is a  sequence $ \{f_{n} \}_{n=1}^{\infty} $ with $ {f_{n}}\to {f}$ locally uniformly and with $\limsup_{n\to\infty} H( {f_{n}})<H( {f})$. Our main result  shows this is  true for affine mappings. Addressing conjectures of F.W. Gehring and Iwaniec we show the jump up in the limit can be arbitrarily large and give conjecturally sharp bounds :  for each $\alpha<\sqrt{2}$ there is  ${f_{n}}\to {f}$ locally uniformly  with $f$ affine and 
\[ \alpha \; \limsup_{n\to\infty} H( {f_{n}}) <   H( {f})  \]
We conjecture $\sqrt{2}$ to be best possible.
}
\end{abstract} 

\section{Introduction.}    
This article is concerned with the convexity properties of the linear distortion functional and in particular its lowersemicontinuity. We therefore begin with a definition. 

\begin{definition} Let $\Omega$ be a domain in $\IR^n$, $n\geq 2$, and $f:\Omega \to f(\Omega)\subset\IR^n$ a homeomorphism.  For each $x\in \Omega$ and $0<r<d(x,\partial\Omega)$  set
\begin{equation}
H(x,f) = \limsup_{r\to 0} \;\;\frac{\max_{|x-y|=r}\{|f(x)-f(y)|\}}{\min_{|x-y|=r}\{|f(x)-f(y)|\}}.
\end{equation}
If $H(x,f)$ is bounded in $\Omega$, and if $H=H(f)=\|H(x,f)\|_{L^\infty(\Omega)}$, then we say $f$ is {\em $H$-quasiconformal} and $H$ is the {\em linear distortion} of $f$.
\end{definition}
If $f$ has a nonsingular derivative $Df(x)$ at $x\in \Omega$ with singular values $\lambda_1(x)\leq \lambda_2(x) \leq \cdots \leq \lambda_n(x)$, then one can see
\[ H(x,f)=\frac{|\lambda_n(x)|}{|\lambda_1(x)|}. \]
It is a remarkable fact that $H(f)$ bounded in $\Omega$ implies Sobolev regularity $f\in W^{1,n}(\Omega)$,  \cite{vaisala,GMP}. Even more remarkable is the result which states $\limsup$ can be replaced by $\liminf$ in the definition when $\Omega=\IR^n$, \cite{HK}.

Other definitions assume this regularity and define the distortion in terms of the differential matrix.  Two common definitions are the inner distortion and  outer distortion,  but there are many others, see \cite[Chapter 9]{IM}.  As examples for a homeomorphism of Sobolev cass $W^{1,n}(\Omega)$ set
\[ K_0(f) = \left\| \frac{|Df|^n}{J(x,f)} \right\|_{L^\infty(\Omega)}, \;  \IK(f) = \left\| \frac{\|Df\|^{n}}{J(x,f)} \right\|_{L^\infty(\Omega)}, \; K_I(f)=\left\|\frac{J(x,f)}{\lambda_1^n}\right\|_{L^\infty(\Omega)} \]
Each distortion defines the same class of maps,  and easy eigenvalue calculations show that
\begin{equation}\label{compare} \IK(f) \leq  K_0(f) \leq H(f)^{1/n}, \quad  H(f) \leq K_I(f) \leq H(f)^n.  \end{equation}
Polyconvexity of the distortion functional (when the distortion is a convex function of minors of the differential) and the usual compactness properties of quasiconformal mappings \cite{GM} show that each of the distortions $\IK(f),K_0(f)$ and $K_I(f)$ has the following lower semicontinuity property.
\begin{theorem}\label{Thm1} Let $K=K(f)$ denote one of the three distortion functionals $\IK(f)$, $K_0(f)$, $K_I(f)$, let $K_\infty<\infty$ and let $\{f_j:\Omega\to \IR^n\}$ be a sequence of quasiconformal mappings with $K(f)\leq K_\infty$.  Then either
\begin{enumerate}
\item The sequence $f_j$ tends locally uniformly in $\Omega$ to a constant mapping with value in $\overline{\IR^n}$.
\item  There is $x_0\in \Omega$ and  $f_j$ tends locally uniformly in $\Omega\setminus \{x_0\}$ to a constant mapping with value in $\overline{\IR^n}$.
\item  The sequence $f_j$ tends locally uniformly in $\Omega$ to a $K(f)$ quasiconformal mapping $f_\infty:\Omega\to\IR^n$ and $K(f)\leq K_\infty$.  
\end{enumerate}
\end{theorem}

Since the very beginning of the multidimensional theory of quasiconformal mappings it was widely believed that the class of $H$-quasiconformal mappings in $\mathbb{R}^n$, defined via the linear distortion,  is closed with respect to local uniform convergence. However this question of lower semicontinuity was answered negatively by Tadeusz Iwaniec in response to a query of Curt McMullen. In his paper \cite{I}, he gave an explicit example to prove the following.
\begin{theorem}
There exists a sequence of $H$-quasiconformal mappings $f_j:\IR^n\to\IR^n$ converging locally uniformly in $\IR^n$ to a linear mapping $f:\IR^n\to\IR^n$ with
\begin{equation} 
H(f) >  \limsup_{\nu\to\infty} H(f_{j}). 
\end{equation}
\end{theorem}
As we will see, the reason for this unusual and anomalous behaviour of the linear distortion function is that it fails to be rank-one convex in dimensions higher than 2. A natural question is how big the jump up can be here and Gehring and Iwaniec give the following bounds found through Theorem \ref{Thm1} and (\ref{compare}).

\begin{theorem}\label{thmgi} Suppose that $f_j: \Omega \to \IR^n$  is a sequence of quasiconformal mappings which converges weakly in $W^{1,n}(\Omega)$ to $f$ and suppose that
$H(x,f_j) \leq  M$ in $\Omega$ for $j= 1,2,\ldots$. Then 
\begin{equation} \label{hbound}
H(x,f)\leq \frac{1}{2} \big(M+M^{n-1}\big)^{2/n}  
\end{equation}
 \end{theorem}
 When $n=3$ and $M$ is large this bound is roughly $\frac{1}{2}M^{4/3}$ while it is $\frac{1}{2}M^2$ for $n$ large. We show the best possible lower bound for the right-hand side of (\ref{hbound}) here must exceed $\sqrt{2} M$ for $M$ large,  with explicit bounds for all $M$.
 
 \section{Main results} Our main results are the following which give lower bounds for Gehring and Iwaniec's result,  and establish the generic nature of the failure of lower semicontinuity for the linear distortion -- at least among affine mappings.
 
 \begin{theorem} Let $A:\IR^n\to\IR^n$ be an affine mapping whose differential has three distinct singular values.  Then there is a sequence of $H$--quasiconformal mappings $f_j:\IR^n\to \IR^n$,  $H(x,f_j)\leq H$, which converge to $A$ uniformly in the spherical metric and $H<H(x,A)=H(A)$.
 \end{theorem}
We give explicit bounds for the jump here in terms of the singular values of $A$, though they are complicated.  The next result relates to Theorem \ref{thmgi}.
 
 \begin{theorem} Let $\alpha<\sqrt{2}$.  Then there is a sequence of $H$-quasiconformal mappings $f_j:\IR^n\to\IR^n$ with $H(x,f_j)\leq H$ converging locally uniformly to an affine mapping and
 \begin{equation}
 H(A)=H(x,A) \geq \alpha H.
 \end{equation}
 \end{theorem}
For our examples $H \to \infty$ as $\alpha\to \sqrt{2}$ and so we see that the gap $H(A)-H(f_j)$ can be arbitrarily large.  We give strong numerical evidence to suggest that $\sqrt{2}$ is optimal here,  at least in three dimensions.    

 \medskip
 
 These two results are based on the following properties of the linear distortion functional.  As an example we 
 recall that the determinant function $\det : \IR^{n\times n} \to \IR$, in spite of the non-linearity of this polynomial of $n^2$ variables, is in fact linear in the directions of rank-one matrices. More precisely,  the function of the real variable $t \mapsto \det (A+tB)$ is linear if $\mathrm{rank}(B)\leq 1$. The same is true for lower-order minors and consequently for null-Lagrangians, being linear combinations of the minors of $Df$. 
 
A rank-one matrix  can be written as the tensor product of two vectors.  The key idea in Iwaniec's work is that the linear distortion function fails to be rank-one convex in dimension $ n \geq 3 $.  To compute the linear distortion of an affine mapping
$x\mapsto Ax+b$ we study the eigenvalues of
$$ A^t A \in Sym^{+}_{3\times 3} (\IR), $$
the space of symmetric positive definite $3\times 3$ matrices. Given such an $A$ the spectral theorem tells us $A$ is orthogonally  diagonalisable.  It is an elementary fact that if $U,V$ are orthogonal and $f$ is quasiconformal,  then $H(x,f)=H(UfV,V^{-1}x)$,  and so we may as well suppose  $A$ is  diagonal.   In this way, we reduce the problem of the convexity of the linear distortion functional to considering that functional defined on the space of $3\times 3 $ diagonal matrices with entries $1=a_{11}\leq a_{22} \leq a_{33}$.  Iwaniec gave an elementary argument to go from three-dimensions to $n$-dimensions which we later  recall.  

\medskip

To achieve these explicit bounds we study an interesting question of independent interest and possibly connected with some aspects of materials science : determine the rank-one direction  for which the linear distortion function at $A$ is ``most concave''. These directions  might identify the structure of the laminations for the minimisers of certain stored energy functionals occurring in the calculus of variations, \cite{MG}. Thus we frame our proof through two problems we address.

\subsection{Problem 1 : Best rank-one direction.} Let $$A=\begin{bmatrix} 1 & 0 & 0 \\ 0 & a & 0 \\ 0 & 0 & b \end{bmatrix}={\rm diag}(1,a,b)$$ be {\rm diag}onal with $1<a<b$.  Determine vectors $\mathbf{u}$ and  $\mathbf{v}$ in $ \IR^{3} $, $\|\mathbf{u}\|=\|\mathbf{v}\| =1$, so that with $B_0=\mathbf{u} \otimes \mathbf{v} $   we have 
\begin{eqnarray}\label{6}
\frac{d}{dt}\Big|_{t=0} H(A+t\, B_{0})  & = &  0, \\
\frac{d^2}{dt^2}\Big|_{t=0} H(A+t\, B_{0})& \leq & \frac{d^2}{dt^2}\Big|_{t=0} H(A+t\, B), \label{7}
\end{eqnarray} 
for every rank-one matrix $B= \tilde{\mathbf{u}}\otimes\tilde{\mathbf{v}}  \in \IR^{3\times 3}$, $\|\tilde{\mathbf{u}}\|=\|\tilde{\mathbf{v}}\| =1$ with 
\[ \frac{d}{dt}\Big|_{t=0} H(A+t\, B)=0.\]

The solution to Problem 1 is unique up to sign.  It is in this direction we might expect to find the minimum values of $H(A+tB)$.
Next we identify the $t$-interval that $H(A+t\, B)$ is concave.
 
\subsection{Problem 2: Intervals of concavity} Let $A={\rm diag}(1,a,b)$  with $1<a<b$ and suppose  $B_0=\mathbf{u} \otimes \mathbf{v}^t $   is a solution to Problem 1.  Determine the largest real numbers $\mathbf{t}_+>0$ and $\mathbf{t}_-<0$ so that $ H(A+t\, B_0) $ is a smooth function of $t$ in the interval $\mathbf{t}_- <t<\mathbf{t}_+$.  Then
determine $H(A+\mathbf{t}_- B_0)$ and $H(A+\mathbf{t}_+ B_0)$. 

\medskip

We naturally expect that the values $\mathbf{t}_-$ and $\mathbf{t}_+$ are where the singular values of $A+t B_0$ cross as $t$ varies - they must cross as $H(A+t B_0)\to +\infty$ as $t\to \pm \infty$. The values $\mathbf{t}_- $   and $\mathbf{t}_+$ will be determined from a (rather challenging) discriminant problem.  We conjecture that
for all rank-one matrices $B=\tilde{u}\otimes\tilde{v}$, $\|\tilde{u}\|=\|\tilde{v}\|=1$ with 
$$ \frac{\rm d}{\rm{d}t}\Big|_{t=0} H(A+t\, B)=  0, \;\;\; {\rm and} \;\;\; \frac{\rm d^2}{\rm dt^2}\Big|_{t=0} H(A+t\, B) < 0, $$ 
we have for all $t>0$
$$ H(A+t B) \geq \max \{ H(A+\mathbf{t}_- B_0),  H(A+\mathbf{t}_+ B_0)\}. $$
This conjecture expresses the hope that the ``best rank-one direction'' also leads to the largest gap between $H(A)$ and $H(A+t\, B)$ and therefore gives us the approximation to $A$ of least linear distortion. However, it may be that there is another path giving a better result for larger $t$, though our numerical evidence suggests otherwise.

\section{Solving Problem 1. } Our problem now is to determine the best rank-one matrix $B_0=\mathbf{u}_{0}.\mathbf{v}_{0}^t$ so that (\ref{6}) and (\ref{7}) hold, That is the coefficient of the quadratic term in the series expansion of $H(A+tB)$ is as negative as possible.  Let $A={\rm diag}(1,a,b)$, $1< a<b$, and $B$ be a rank-one matrix, $ B=\mathbf{u} \otimes \mathbf{v}=\mathbf{u}\cdot \mathbf{v}^{t} $, where 
\[ \mathbf{u}^t= (\sqrt{1-r^2},r\cos(\theta_1),r\sin(\theta_1)), \quad \mathbf{v}^t= (\sqrt{1-s^2},s\cos(\theta_2),s\sin(\theta_2))\] 
Here $0\leq r,s \leq 1$ and  $ \theta_{1},\theta_{2}\in[0,2\pi]$.  Clearly $H(A)= b$. Let $\lambda_{1}(t)$, $\lambda_{2}(t)$ and $\lambda_{3}(t)$ be eigenvalues of $X=(A+tB)^t(A+tB)$. The functions $\lambda_{i}(t)$ are locally well defined and smooth in $t$ and $\lambda_{1}(0)= 1,  \lambda_{2}(0)= a^2,  \lambda_{3}(0)= b^2.$
Thus for sufficiently small $t$  we have
$\lambda_{1}(t) < \lambda_{2}(t) < \lambda_{3}(t) $ and
$$ H(A+tB)=\sqrt{\frac{\lambda_{3}(t)}{\lambda_{1}(t)}}. $$
As $$(\frac{A}{t}+B)^t (\frac{A}{t}+B)=\frac{1}{t^2}A^tA+\frac{1}{t}(A^tB+B^tA)+B^tB, $$ if $t$ tends to infinity then 
$$ \lim_{|t| \rightarrow +\infty} H(A+tB)=\lim_{|t| \rightarrow +\infty} H(\frac{A}{t}+B)=H(B)= \infty,$$
as $B$ is  rank-one.   
 The matrix $B$ has the form
$$ B=\begin{bmatrix}
 \sqrt{1-r^2} \sqrt{1-s^2} & \sqrt{1-r^2} s \cos (\theta_2) & \sqrt{1-r^2} s \sin (\theta_2) \\
 r \sqrt{1-s^2} \cos (\theta_1) & r s \cos (\theta_1) \cos (\theta_2) & r s \cos (\theta_1) \sin (\theta_2) \\
 r \sqrt{1-s^2} \sin (\theta_1) & r s \cos (\theta_2) \sin (\theta_1) & r s \sin (\theta_1) \sin (\theta_2) 
 \end{bmatrix} $$
and $A+tB$ has the form
 $$
\begin{bmatrix}
  1+t\sqrt{1-r^2} \sqrt{1-s^2} & ts\sqrt{1-r^2}  \cos (\theta_2) &ts \sqrt{1-r^2}  \sin (\theta_2) \\
t r \sqrt{1-s^2} \cos (\theta_1) &a+ t r s \cos (\theta_1) \cos (\theta_2) &t r s \cos (\theta_1) \sin (\theta_2) \\
t r \sqrt{1-s^2} \sin (\theta_1) &t r s \cos (\theta_2) \sin (\theta_1) &b+t r s \sin (\theta_1) \sin (\theta_2) 
 \end{bmatrix}.$$  
Using a second order Taylor series in $t$ we may find the smallest and the largest eigenvalues of $X$ to second order. Let $I$ be the identity $3\times 3$ matrix. Then the smallest eigenvalues of $X$ can be found from 
\begin{equation}\label{2ndorder} \mathrm{det}[X-\lambda_{1} I] \approx \mathrm{det}[A+tB-(1+xt+yt^{2})I]=0, \end{equation}
where $ \lambda_{1}<\lambda_{2}<\lambda_{3}$. Differentiating (\ref{2ndorder}) gives $x= 2 \sqrt{1-r^2} \sqrt{1-s^2}$ and
\begin{eqnarray*}
y&=& 1-s^2-\frac{\Big(ar\sqrt{1-s^2}\cos(\theta_1)+s\sqrt{1-r^2}\cos(\theta_2) \Big)^2 }{a^2-1} \\&& - \frac{\Big(br\sqrt{1-s^2}\sin(\theta_1) s\sqrt{1-r^2}\sin(\theta_2) \Big)^2 }{b^2-1}.
\end{eqnarray*}
Hence
\begin{footnotesize}
\begin{equation*}
\begin{aligned} 
\lambda_{1}({t})=&1 +2 \big(\sqrt{1-r^2} \sqrt{1-s^2}\big)  {t} +\Big(1-s^2-\frac{\big(ar\sqrt{1-s^2}\cos(\theta_1)+s\sqrt{1-r^2}\cos(\theta_2) \big)^2 }{a^2-1} - \\&
\frac{\big(br\sqrt{1-s^2}\sin(\theta_1)+s\sqrt{1-r^2}\sin(\theta_2) \big)^2 }{b^2-1} \Big){t^2} +O( {t^3}) . 
\end{aligned}
\end{equation*}
\end{footnotesize}
Similarly the largest eigenvalue to second order is
\begin{small}
\begin{equation*}
\begin{aligned} 
\lambda_{3}({t})&=b^2 + 2\;brs\sin(\theta_1)\sin(\theta_2)\;  {t}+
 \Big(\frac{s^2}{2}-\frac{1}{2}s^2 \cos(2\theta_2)+\\&\frac{\big(br\sqrt{1-s^2}\sin(\theta_1)+s\sqrt{1-r^2}\sin(\theta_2) \big)^2 }{b^2-1}
 +\frac{s^2\big(rb \:\cos(\theta_2)\sin(\theta_1)+ra \:\cos(\theta_1)\sin(\theta_2) \big)^2}{b^2 -a^2}   \Big)\; {t^2}.
\end{aligned}
\end{equation*}
\end{small}
Therefore for small enough $t$, the linear distortion function $H(A+{t}B)$ is    
\begin{small} 
\begin{eqnarray*}
&&  b+\big(-b\sqrt{1-r^2} \sqrt{1-s^2}+rs \sin(\theta_1)\sin(\theta_2) \big) {t}+\frac{1}{2b}\Big(\frac{s^2}{2}-\frac{1}{2}s^2 \cos(2\theta_2) \\
&&
- 4brs\sqrt{1-r^2} \sqrt{1-s^2}\:\sin(\theta_1)\sin(\theta_2)+\frac{\big(br\sqrt{1-s^2}\sin(\theta_1)+s\sqrt{1-r^2}\sin(\theta_2) \big)^2 }{b^2-1}\\
&&+
\frac{s^2\big(rb \:\cos(\theta_2)\sin(\theta_1)+ra \:\cos(\theta_1)\sin(\theta_2) \big)^2}{b^2 -a^2}-\big(b\sqrt{1-r^2} \sqrt{1-s^2}-rs \sin(\theta_1)\sin(\theta_2) \big)^2\\ 
&&+ 
b^2\Big(-1 +s^2+4(r^2-1)(s^2-1)+\frac{\Big(ar\sqrt{1-s^2}\cos(\theta_1)+s\sqrt{1-r^2}\cos(\theta_2) \Big)^2 }{a^2-1}\\
&& +
\frac{\big(br\sqrt{1-s^2}\sin(\theta_1)+s\sqrt{1-r^2}\sin(\theta_2) \big)^2 }{b^2-1} \Big) \Big) {t^2}+O( {t^3}).
\end{eqnarray*}
\end{small} 
We require the first derivative of $H(A+tB)$ to be zero.  
\begin{equation}\label{rs} b\big(\sqrt{1-r^2} \sqrt{1-s^2}\:\big)=rs \sin(\theta_1)\sin(\theta_2). \end{equation}
We want to minimise the quadratic coefficient above as a function of the four variables $r$, $s$, $\theta_{1}$ and $\theta_{2}$ and the two parameters $a$ and $b$. $Q(r,s,\theta_1,\theta_2)=$
\begin{small}
\begin{eqnarray*}
&& \frac{1}{2b}\Big(\frac{s^2}{2}-\frac{1}{2}s^2 \cos(2\theta_2)- 4brs\sqrt{1-r^2} \sqrt{1-s^2}\:\sin(\theta_1)\sin(\theta_2)\\ 
& &+
\frac{\big(br\sqrt{1-s^2}\sin(\theta_1)+s\sqrt{1-r^2}\sin(\theta_2) \big)^2 }{b^2-1}+\frac{s^2\big(rb \:\cos(\theta_2)\sin(\theta_1)+ra \:\cos(\theta_1)\sin(\theta_2) \big)^2}{b^2 -a^2} \\ & &-
\big(b\sqrt{1-r^2} \sqrt{1-s^2}-rs \sin(\theta_1)\sin(\theta_2) \big)^2 +b^2\Big(-1 +s^2+4(r^2-1)(s^2-1)\\ & &+
\frac{\Big(ar\sqrt{1-s^2}\cos(\theta_1)+s\sqrt{1-r^2}\cos(\theta_2) \Big)^2 }{a^2-1}+\frac{\big(br\sqrt{1-s^2}\sin(\theta_1)+s\sqrt{1-r^2}\sin(\theta_2) \big)^2 }{b^2-1} \Big) \Big). 
\end{eqnarray*}\end{small}
We substitute $r^2s^2$ using (\ref{rs}) to find $Q(r,s,\theta_1,\theta_2)=$
\begin{small}
\begin{eqnarray*}
& & \frac{1}{8} \Big(\frac{r^2\cos(2\theta_1)\big(-2(b^2+b^4)+s^2(-5+3b^2+2b^4)-(b^2-5)s^2\cos(2\theta_2) \big)}{b^3-b}\\&&+
 \frac{3r^2 s^2+2b^4(5r^2-4)(s^2-1)+b^2\big(r^2(14-17s^2)+12s^2-8\big)}{b^3-b}\\&&+\frac{(-4b^2+3(b^2-1)r^2 )s^2 \cos(2\theta_2)}{b^3-b}\\&&+ 
  4bs^2\cos^2(\theta_2)\Big(\frac{1-r^2}{a^2-1}+\frac{r^2 \sin^2(\theta_1)}{b^2-a^2} \Big)+4a^2 r^2 \cos^2(\theta_1)\Big(\frac{b(1-s^2)}{a^2-1}+\frac{s^2 \sin^2(\theta_1)}{b^3-a^2 b} \Big)\\&&-
 \frac{2a(b^2-1)r^2 s^2\sin (2 \theta_1) \sin (2 \theta_2)}{(a^2-1)(a^2-b^2)}  \Big).
\end{eqnarray*} 
\end{small}
It is obvious that the function $Q$ is $\pi$-periodic, so we assume the values of $\theta_1,\theta_2\in[0 \pi]$. Put {$\delta=r^2$} and {$\eta=s^2$} and write the equation as 
\begin{small}\[
\begin{aligned}
&Q({\delta},{\eta},\theta_1,\theta_2)=\frac{1}{8} \Big(\frac{{\delta}\cos(2\theta_1)\big(-2(b^2+b^4)+{\eta}(-5+3b^2+2b^4)-(b^2-5){\eta}\cos(2\theta_2) \big)}{b^3-b}+ \\& 
 \frac{3{\delta} {\eta}+2b^4(5{\delta}-4)({\eta}-1)+b^2\big({\delta}(14-17{\eta})+12{\eta}-8\big)+(-4b^2+3(b^2-1){\delta} ){\eta} \cos(2\theta_2)}{b^3-b}+\\& 
  4b{\eta}\cos^2(\theta_2)\Big(\frac{1-{\delta}}{a^2-1}+\frac{{\delta} \sin^2(\theta_1)}{b^2-a^2} \Big)+4a^2 {\delta} \cos^2(\theta_1)\Big(\frac{b(1-{\eta})}{a^2-1}+\frac{{\eta} \sin^2(\theta_1)}{b^3-a^2 b} \Big)-\\& 
 \frac{2a(b^2-1){\delta} {\eta}\sin (2 \theta_1) \sin (2 \theta_2)}{(a^2-1)(a^2-b^2)}  \Big).
\end{aligned}\]
\end{small}
Now 
\[ \delta \eta=\frac{b^2\big(\delta+\eta -1 \big)}{b^2- \sin^2(\theta_1)\sin^2(\theta_2)} \] by (\ref{rs}) so we can eliminate the nonlinear term $\delta\eta$ and the function $Q$ can be simplified to $Q(\delta,\eta,\theta_1,\theta_2)=$
\begin{scriptsize} 
\begin{align*}
& -\frac{1}{32 \left(a^2-1\right)\left(b^2-1\right) \left(b^2-a^2\right) (b^2-\sin^2(\theta _1) \sin ^2(\theta _2))}\\ & 
  b \Big(b^4(32-7\delta-7\eta)+8b^6(\delta+\eta-2)-7b^2(\delta+\eta)+2a^4\big(\delta+\eta-4-4b^2(\delta+\eta-3)\big)\\ &-
a^2\big(\delta-8+b^2(40-21(\delta+\eta))+\eta+8b^4(\delta+\eta)\big)+(a^2-b^2)\Big[\Big(\eta-8-2a^2(-4+4b^2(1+\delta-\eta)+\eta)\\ &+
b^2(8\delta-8b^2(\eta-1)+\eta)\Big)\cos(2\theta_1)+(1-2a^2+b^2)\delta \cos(4\theta_1) \Big]+\Big[8(a^2-1)\big(-b^2(\delta+\eta-2)\\ &+
a^2(-1+b^2(\delta+\eta-1))\big)\cos(2\theta_1)+(a^2-b^2)\Big(-8+\delta+a^2(8-2\delta+8b^2(\delta-\eta-1)+b^2(-8b^2(\delta-1)\\ &+
+\delta+8\eta)+(2a^2-b^2-1)\delta\cos(4\theta_1)\Big)\Big]\cos(2\theta_2)-2(a^2-b^2)(2a^2-b^2-1)\eta\cos(4\theta_2)\sin^2(\theta_1)\\ &-                   
  8ab(b^2-1)^2(\delta+\eta-1)\sin(2\theta_1)\sin(2\theta_2) \Big).
 \end{align*}
\end{scriptsize}
$\xi=32(a^2-1)(b^2-1)(b^2-a^2)$ and $\mu=\big(b^2-\sin^2(\theta_1)\sin^2(\theta_2)\big)$ are positive. 
We rewrite $Q$ with respect to four variables $\delta$, $\eta$, $\theta_1$ and $\theta_2$ as below, 
\begin{equation}\label{equation 4.15}
Q({\delta},{\eta},\theta_1,\theta_2)=\frac{ {\alpha}{\delta}+{\beta}{\eta}+{\gamma}}{32(a^2-1)(b^2-1)(b^2-a^2)\big(b^2-\sin^2(\theta_1)\sin^2(\theta_2)\big)}.\vspace{0.4 cm}
\end{equation}
 where
\[
\begin{aligned}
{\alpha(\theta_1,\theta_2)}=& b \Big(a^2-2a^4+\big(7-21a^2+8a^4\big)b^2+\big(7+8a^2\big)b^4-8b^6- 8\big(a^2-1\big)\\&
b^2\cos(2\theta_1)\Big(-a^2+b^2+\big(a^2-1\big)\cos(2\theta_2)\Big)-\big(a^2-b^2\big)\Big[\Big(1+b^2-8b^4  \\&
+a^2\big(-2+8b^2\big)\Big)\cos(2\theta_2)+2\Big(1-2a^2+b^2\Big)\cos(4\theta_1)\sin^2(\theta_2)\Big]+\\&
8ab\Big(b^2-1\Big)^2\sin(2\theta_1)\sin(2\theta_2)\Big),       
\end{aligned}
\]
\[
\begin{aligned}
{\beta(\theta_1,\theta_2)}=& b \Big(a^2-2a^4+\big(7-21a^2+8a^4\big)b^2+\big(7+8a^2\big)b^4-8b^6+\cos(2\theta_1)\\&
\Big[-\big(a^2-b^2\big)\Big(1+b^2-8b^4+a^2\big(8b^2-2\big)\Big)-8\big(a^2-1\big)b^2\cos(2\theta_2)\Big]+\\&
2\big(b^2-a^2\big)\Big[-4\big(a^2-1\big)b^2\cos(2\theta_2)+\Big(1-2a^2+b^2\Big)\cos(4\theta_2)\sin^2(\theta_1)\Big]\\&
+8ab\Big(b^2-1\Big)^2\sin(2\theta_1)\sin(2\theta_2)\Big),       
\end{aligned}]
\]
\[
\begin{aligned}
{\gamma(\theta_1,\theta_2)}=& -b \Big(32b^4 -16b^6 +8a^2\big(1-5b^2\big)+8a^4\big(3b^2-1\big)-8\big(a^2-b^2\big)\\&
\big(a^2-b^2-1\big)\cos(2\theta_1)+\cos(2\theta_2)\Big[8\big(b^2-1\big)\big(a^2-b^2\big)\big(1-a^2+b^2\big)-\\&
8\big(a^2-1\big)\Big(a^2+\big(a^2-2\big)b^2\Big)\cos(2\theta_1)\Big]+8ab\Big(b^2-1\Big)^2\sin(2\theta_1)\sin(2\theta_2)\Big).       
\end{aligned}\]
The constraint is given at (\ref{rs}), and since the denominator of $Q$ does not vanish, we may multiply the constraint by this term and clear the multiplicative factor. We then consider Lagrange multipliers to examine the following function.  
 \vspace{0.5cm}
\begin{equation}\label{equation 4.20}
F({\delta},{\eta},\lambda)=\frac{ {\alpha}{\delta}+{\beta}{\eta}+{\gamma}}{\xi\mu}-\lambda \big(b^2(1-{\delta}-{\eta})+{\delta} {\eta} \mu \big). 
\end{equation}
We get the following three equations  
\begin{enumerate}
\item $\frac{\partial F}{\partial \delta}=b^2 \lambda -\eta \lambda \mu +\frac{\alpha}{\xi\mu}=0,$
\item $\frac{\partial F}{\partial \eta}=b^2 \lambda -\delta \lambda \mu +\frac{\beta}{\xi\mu}=0,$
\item $\frac{\partial F}{\partial \lambda}=b^2(-1+\delta+\eta)-\delta \eta\mu=0.$
\end{enumerate}
giving the two sets of solutions for $\delta$ and $\eta$. 
\[
\delta_1=\frac{b\Big(b-\sqrt{\frac{\beta}{\alpha}}\sqrt{b^2-\mu}\Big)}{\mu},\; \eta_1=\frac{b\Big(b-\sqrt{\frac{\alpha}{\beta}}\sqrt{b^2-\mu}\Big)}{\mu},\; \lambda_1=-\frac{\sqrt{\alpha\beta}}{b\:\xi\:\mu\sqrt{b^2-\mu}},
\]
\[
\delta_2=\frac{b\Big(b+\sqrt{\frac{\beta}{\alpha}}\sqrt{b^2-\mu}\Big)}{\mu},\; \eta_2=\frac{b\Big(b+\sqrt{\frac{\alpha}{\beta}}\sqrt{b^2-\mu}\Big)}{\mu},\;\lambda_2=\frac{\sqrt{\alpha\beta}}{b\:\xi\:\mu\sqrt{b^2-\mu}}.
\]
Since $0 \leq \delta, \:\eta \leq 1$, we must check which is the set of solutions between $0$ and $1$ that we want. In fact $0 \leq \delta_1 \leq 1$,  $0 \leq \eta_1 \leq 1$, $ \delta_2 \geq 1$ and $\eta_2 \geq1$. We have 
$\mu=\big(b^2-\sin^2(\theta_1)\sin^2(\theta_2)\big)$, so  
$$\sqrt{b^2-\mu}=\sqrt{\sin^2(\theta_1)\sin^2(\theta_2)}=\sin(\theta_1)\sin(\theta_2) \geq 0.$$
For $\delta_2$, we have  
$$ b^2+b\sqrt{\frac{\beta}{\alpha}}\:(\sin(\theta_1)\sin(\theta_2))  \geq b^2 \geq  b^2-\sin^2(\theta_1)\sin^2(\theta_2) . $$
So, 
$$\delta_2=\frac{ b^2+b\sqrt{\frac{\beta}{\alpha}}\:(\sin(\theta_1)\sin(\theta_2)) }{ b^2-\sin^2(\theta_1)\sin^2(\theta_2) }\geq 1. $$ 
Similarly
$$\eta_2=\frac{ b^2+b\sqrt{\frac{\alpha}{\beta}}\:(\sin(\theta_1)\sin(\theta_2)) }{ b^2-\sin^2(\theta_1)\sin^2(\theta_2) }\geq 1.  $$ 
\subsubsection{An extremal case.}\label{goodcase} We next claim the minimum of the function $Q$ is negative by examining a special case which will turn out to be the extremal.  
Let $r=s$ and $\theta_2=\pi-\theta_1$. Then the vectors $\mathbf{u}$ and $\mathbf{v}$ are 
$$\mathbf{u}^t= (\sqrt{1-r^2},r\cos(\theta_1),r\sin(\theta_1)),\quad \mathbf{v}^t= (\sqrt{1-r^2},-r\cos(\theta_1),r\sin(\theta_1)).$$
Thus 
$$ A+tB=\begin{bmatrix}
1+t(1-r^2) & -t r\sqrt{1-r^2}  \cos (\theta_1) &tr \sqrt{1-r^2}  \sin (\theta_1) \\
t r \sqrt{1-r^2} \cos (\theta_1) &a- t r^2 \cos^2 (\theta_1) &t r^2 \cos (\theta_1) \sin (\theta_1) \\
t r \sqrt{1-r^2} \sin (\theta_1) &-t r^2 \cos (\theta_1) \sin (\theta_1) &b+t r^2 \sin^2 (\theta_1) 
 \end{bmatrix}. $$  
The smallest and the largest eigenvalues of the matrix $(A+tB)^t(A+tB)$ can be found to second order as before.
\begin{small} 
\begin{equation*}
\begin{aligned} 
\lambda_{1}(t)=&1 +2 (1-r^2)  t +\Big(\frac{(r^2-1)\big(1+a-b+r^2+ab\,(r^2-1)-(a+b)r^2\cos(2\theta_1) \big)}{(a+1)(b-1)}  \Big)t^2, 
\end{aligned}
\end{equation*}
\end{small} 
and 
\begin{small}
\begin{equation*}
\begin{aligned} 
\lambda_{3}(t)&=b^2 + 2\;br^2\sin^2(\theta_1)\;  t-\\
&\Big(\frac{r^2\Big(b\big(3r^2+b(r^2-4)\big)+a\big(r^2+b(3r^2-4)\big)+(a-b)(b-1)\cos(2\theta_1)\Big)\sin^2(\theta_1)}{2(a+b)(b-1)}  \Big)t^2.
\end{aligned}
\end{equation*}
\end{small}
The linear distortion function will be 
$$ H(A+tB)=\sqrt{\frac{\lambda_{3}(t)}{\lambda_{1}(t)}}=b+L(r,\theta_1)\,t+Q(r,\theta_1)\,t^2+ O(t^3).$$ 
The first derivative of the linear distortion function must be zero,
$L(r,\theta_1)= b(r^2-1) +r^2 \sin^2(\theta_1)=0,$
so  
$b(r^2-1)=- r^2 \sin^2(\theta_1)$.
If we put the above equation in the quadratic coefficient function we see $Q(r, \theta_1)$ is
\begin{equation*}
\begin{aligned} &\frac{r^2\Big(4(a+1)(a+b)-\big(1+a+a^2+(a-1)b+b^2  \big)r^2  \Big)}{8(a+1)(b-1)(a+b)}\\
&+\frac{r^2\Big(-4(a+1)(a+b)\cos(2\theta_1)+\big(1+a+a^2+(a-1)b+b^2  \big)r^2 \cos(4\theta_1) \Big)}{8(a+1)(b-1)(a+b)}.
\end{aligned}
\end{equation*}
We first find the critical points, $\frac{\partial Q}{\partial r}(r ,\theta_1)=0$  and  $\frac{\partial Q}{\partial \theta_1}(r ,\theta_1)=0$.
The solutions  are $ r\in \{0,1\}$,  $\theta_1\in\{0, \pi\}$.  If $r=0$, then $Q=0$. With $r=1$ the function $Q$ is 
\begin{equation*}
\begin{aligned} 
Q(1, \theta_1)&=\frac{\Big(4(a+1)(a+b)-\big(1+a+a^2+(a-1)b+b^2  \big)  \Big)}{8(a+1)(b-1)(a+b)}\\
&+\frac{\Big(-4(a+1)(a+b)\cos(2\theta_1)+\big(1+a+a^2+(a-1)b+b^2  \big) \cos(4\theta_1) \Big)}{8(a+1)(b-1)(a+b)}.
\end{aligned}
\end{equation*}
In this case the partial derivative of $Q$ with respect to $\theta_1$ is
$$\frac{\partial Q}{\partial \theta_1}(1 ,\theta_1)=\frac{8(a+1)(a+b)\sin(2\theta_1)-4\big(1+a+a^2+(a-1)b+b^2  \big) \sin(4\theta_1) }{8(a+1)(b-1)(a+b)}=0.$$
As $8 (a+b)(b-1)(a+1)>0$ we find five critical points in $[0,\pi]$. 
\[  \theta_1=0,\quad \theta_1=\frac{\pi}{2},\quad \theta_1=\pi, \quad \theta_1=\arctan\frac{y_1}{x_1},\quad \theta_1=\pi-\arctan\frac{y_1}{x_1}, \]
where 
$$x_1=\frac{(a+1)(a+b)}{1+a+a^2-b+ab+b^2},$$
and
$$y_1=\frac{(b-1)\sqrt{1+2a+2a^2+2ab+b^2}}{\sqrt{1+2a+3a^2+a^4-2b+2a^3b+3b^2+3a^2b^2-2b^3+2ab^3+b^4}}.$$
The values of function $Q_1$ at the critical points are 
$$Q_1(0)=0,\qquad\qquad Q_1(\frac{\pi}{2})=\frac{1}{b-1}>0,\qquad\qquad Q_1(\pi)=0,$$
$$ Q_1(\arctan\frac{y_1}{x_1})=-\frac{(b-1)^3}{4(a+1)(b+a)\big(1+a+a^2+b(a-1)+b^2\big)}<0.$$
Now we prove that the quadratic function $Q$ is positive on the boundary where $\theta_1,\theta_2\in\{0,\pi\}$.   So,
$$\theta_1=0\:\; \mathrm{and}\:\; \theta_2=\theta, \qquad\qquad\qquad \theta_1=\pi\:\; \mathrm{and}\:\; \theta_2=\theta,$$
$$\theta_1=\theta\:\; \mathrm{and}\:\; \theta_2=0, \qquad\qquad\qquad \theta_1=\theta\:\; \mathrm{and}\:\; \theta_2=\pi.$$
There are four cases here, and we will prove only the first case  as they are entirely similar. With $\theta_1=0$ and $\theta_2=\theta$  
$$\mathbf{u}^t= (\sqrt{1-r^2},r,0), \quad \mathbf{v}^t= (\sqrt{1-s^2},s\cos(\theta),s\sin(\theta)).$$ As
$-b \sqrt{1-r^2}\sqrt{1-s^2}=0$ either $r=1$ or $s=1$.  If $r=1$, then 
\begin{equation*}
\begin{aligned} 
Q(1,s, \theta)&=\frac{b\big(s^2+a^2(2-3s^2)+2b^2(s^2-1)+(a^2-1)s^2(2\cos^2(\theta) -1)    \big)}{4(a^2-1)(a^2-b^2)}\\
&=\frac{b\big(2s^2(a^2-1)(1-\cos^2(\theta))+2(b^2-a^2)(1-s^2)    \big)}{4(a^2-1)(b^2-a^2)}.
\end{aligned}
\end{equation*}
Since $1<a<b$, $0 \leq s \leq 1$, and $0 \leq \theta \leq \pi$, then   
$ Q(1,s, \theta) \geq 0$. The case $s=1$ is entirely similar.

\medskip
 
We now have shown that $Q$,  the quadratic term, is positive on its boundary with an absolute minimum at $\delta_1$ and $\eta_1$. Note that the equations defining $ \delta_1 $ and $ \eta_1 $ include the expressions 
 $\sqrt{\frac{\beta}{\alpha}}\:\:\: {\rm and} \:\:\:\sqrt{\frac{\alpha}{\beta}},$
 respectively.  Illustrated by Figure 1,  there are computational issues associated with the choice of square roots coming from the way we present the algebraic solution in our formula (which we resolved to present the graphs of the function $Q$ below) so we must be a little careful to account of the ranges of $\beta$ and $ \alpha$ in our calculations.     

\begin{center}
         \includegraphics[scale=0.3]{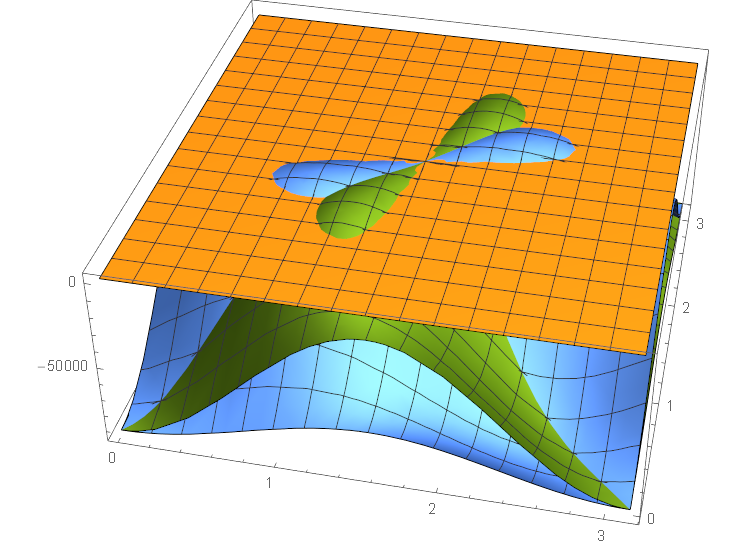} \\
         \end{center}
        {\bf Figure 1.} The functions $\alpha$ and $\beta$ if $A=(1,2,4)$. The functions $\alpha$ and $\beta$ are the blue and the green graphs, respectively.  

Substituting our expressions for $ \delta_1 $ and $ \eta_1 $ and simplifying we find 
$$
Q(\theta_1,\theta_2)=\frac{ \alpha \Big(b^2-b\sqrt{\frac{\beta}{\alpha}}\sin(\theta_1)\sin(\theta_2)\Big)+\beta \Big(b^2-b\sqrt{\frac{\alpha}{\beta}}\sin(\theta_1)\sin(\theta_2)\Big)+\gamma \mu}{\xi\mu^2}.
 $$ 
\begin{center}
      \includegraphics[scale=0.4]{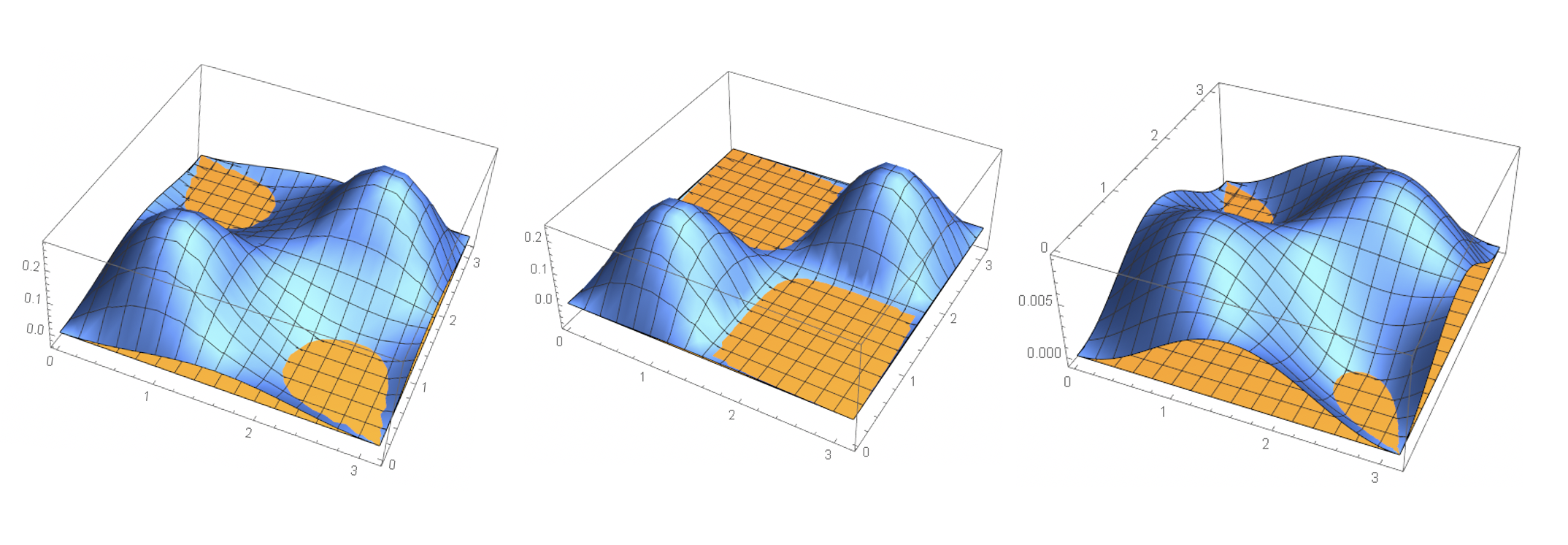} \end{center}
      
      \noindent  {\bf Figure 2} The quadratic function $Q$ if $A=(1,2,10), A=(1,2,105), A=(1,99,154)$ (left to right)

\medskip

The function $Q$ has various properties, for instance: \\
 $Q(\theta_1,\theta_2)=Q(\theta_2,\theta_1),$ 
 $Q(\theta_1,\theta_2)=Q(-\theta_1,-\theta_2),$ 
 $Q(\theta_1,\theta_2)=Q(\pi-\theta_1,\pi-\theta_2),$ 
 $Q(\theta_1,\theta_2)=Q(2\pi-\theta_2,2\pi-\theta_1),$ 
and $Q(\theta_1,\theta_2)=Q(\pi-\theta_2,\pi-\theta_1).$

Some of these are obvious and so we only prove  $Q(\theta_1,\theta_2)=Q(\pi-\theta_1,\pi-\theta_2)$ establishing the symmetry for the functions $\alpha$, $\beta$, $\gamma$ and $\mu$.  
\[
\begin{aligned}
&{\alpha}(\pi-\theta_1,\pi-\theta_2)\\=& b \Big(a^2-2a^4+\big(7-21a^2+8a^4\big)b^2+\big(7+8a^2\big)b^4-8b^6- \\&
8\big(a^2-1\big)b^2\cos(2(\pi-\theta_1))\Big(-a^2+b^2+\big(a^2-1\big)\cos(2(\pi-\theta_2))\Big)-\\&
\big(a^2-b^2\big)\Big[\Big(1+b^2-8b^4  +a^2\big(-2+8b^2\big)\Big)\cos(2(\pi-\theta_2))+\\&
2\Big(1-2a^2+b^2\Big)\cos(4(\pi-\theta_1))\sin^2(\pi-\theta_2)\Big]+\\&
8ab\Big(b^2-1\Big)^2\sin(2(\pi-\theta_1))\sin(2(\pi-\theta_2))\Big)\\& 
=b \Big(a^2-2a^4+\big(7-21a^2+8a^4\big)b^2+\big(7+8a^2\big)b^4-8b^6- 8\big(a^2-1\big)\\&
b^2\cos(2\theta_1)\Big(-a^2+b^2+\big(a^2-1\big)\cos(2\theta_2)\Big)-\big(a^2-b^2\big)\Big[\Big(1+b^2-8b^4  \\&
+a^2\big(-2+8b^2\big)\Big)\cos(2\theta_2)+2\Big(1-2a^2+b^2\Big)\cos(4\theta_1)\sin^2(\theta_2)\Big]+\\&
8ab\Big(b^2-1\Big)^2\sin(2\theta_1)\sin(2\theta_2)\Big) ={\alpha}(\theta_1,\theta_2).    
\end{aligned}
\]
Similarly we see 
\begin{equation}\label{equation 4.26}
\begin{aligned}
&{\beta}(\pi-\theta_1,\pi-\theta_2)={\beta}(\theta_1,\theta_2),\\[0.1 cm]&
{\gamma}(\pi-\theta_1,\pi-\theta_2)={\gamma}(\theta_1,\theta_2).
\end{aligned}
\end{equation}
Also, 
\begin{equation}\label{equation 4.27}
\begin{aligned}
\mu (\pi-\theta_1,\pi-\theta_2)&=\big(b^2-\sin^2(\pi-\theta_1)\sin^2(\pi-\theta_2)\big)\\&
=\big(b^2-(-\sin(\theta_1))^2(-\sin(\theta_2))^2\big)\\&
=\big(b^2-\sin^2(\theta_1)\sin^2(\theta_2)\big)=\mu (\theta_1,\theta_2)
\end{aligned}
\end{equation} 
It now follows that $Q$ has the desired symmetries.   

\medskip

Next, the function $Q$ is positive on the boundary of our region because 
$$ Q(0,\theta_2)=Q(\pi,\theta_2)=\frac{b\sin^{2}(\theta_2)}{2(b^2 - a^2)},\qquad Q(\theta_1,0)=Q(\theta_1,\pi)=\frac{b\sin^{2}(\theta_1)}{2(b^2 - a^2)}. $$ 
To determine the extremes of  $Q$ we must solve
\begin{equation}\label{equation 4.34}
\frac{\partial Q}{\partial \theta_1}(\theta_1,\theta_2)=0, \qquad\qquad \frac{\partial Q}{\partial \theta_2}(\theta_1,\theta_2)=0.
\end{equation}
However these equations are very difficult to solve analytically because they are of eighth degree and trigonometric. Fortunately, it can be proved that the extremes of the function $Q$ are on the lines of symmetry. For this purpose, we must prove that the contour maps of the partial derivatives of $Q$ intersect each other on lines of symmetry, illustrated below in Figure 3. 

 \begin{center}    \includegraphics[scale=0.25]{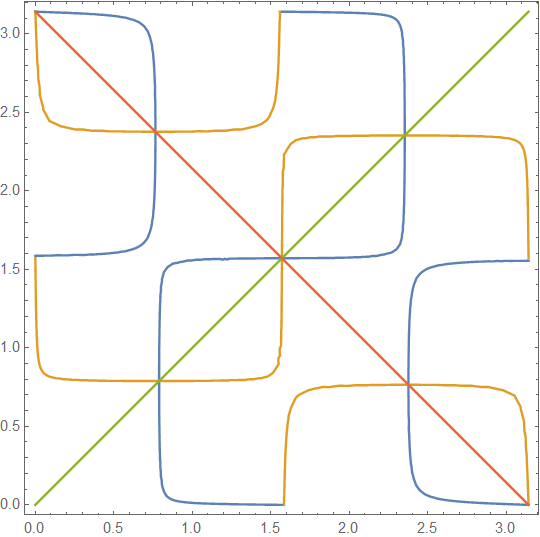} \end{center}
\noindent{\bf Figure 3} The contour maps of $\frac{\partial Q}{\partial \theta_1}(\theta_1,\theta_2)$ and $\frac{\partial Q}{\partial \theta_2}(\theta_1,\theta_2)$, if $A=(1,2,86)$.  
\label{fig: 4.14and4.15}

First, the partial derivatives of $Q$ are symmetric with respect to the line $\theta_1=\theta_2$.  \begin{equation}\label{equation 4.35}
\frac{\partial Q}{\partial \theta_1}(\theta_1,\theta_2)= \frac{\partial Q}{\partial \theta_2}(\theta_2,\theta_1). 
\end{equation}
The point $(\frac{\pi}{2}, \frac{\pi}{2})$ is the intersection point of the lines of symmetry $\theta_1=\theta_2$ and $\theta_1+\theta_2=\pi$. If the origin $(0,0)$ is transferred to the point $(\frac{\pi}{2}, \frac{\pi}{2})$, then the two partial derivatives of the function $Q$ are symmetric with respect to the line $\theta_1=\theta_2$. That is one of these two functions can be considered as $f$ and the other as $f^{-1}$. As a  consequence, the solutions of the simultaneous equations on $[0,\pi]\times [0,\pi]$ are located on the lines of symmetry $\theta_1=\theta_2$ and $\theta_1+\theta_2=\pi$. Hence, the critical points of $Q$ lie on lines of symmetry as illustrated in Figure 3.   
\begin{figure}[ht]
\begin{multicols}{2} 
      \includegraphics[scale=0.385]{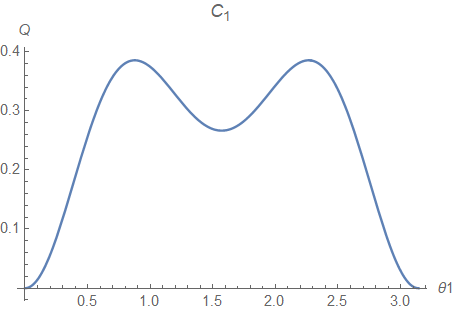}\par
     {{\bf Figure 4.} The cross-section of $Q$ with respect to line $\theta_1=\theta_2$.} 
      \includegraphics[scale=0.34]{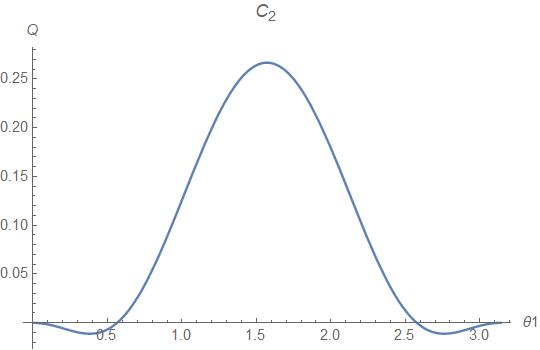}\par
 {{\bf Figure 5.} The cross-section of $Q$ with respect to line $\theta_1+\theta_2=\pi$.}
      \label{fig: 4.17and4.18} 
\end{multicols}
\end{figure}

\subsubsection{Case 1: ($\theta_1=\theta_2$)}  In this case the functions $\alpha$, $\beta$, $\gamma$ and $\mu$ are as follows:
\[
\begin{aligned}
{\alpha(\theta_1)}=&b\Big(a^2-2a^4+\big(7-21a^2+8a^4\big)b^2+\big(7+8a^2\big)b^4-8b^6+\cos(2\theta_1)(a^2-b^2)\\&
(-1+2a^2-9b^2+8b^4)-8(a^2-1)^2 b^2\cos^2(2\theta_1)+2(a^2-b^2)(2a^2-b^2-1) \\&
\cos(4\theta_1)\sin^2(\theta_1)+8ab(b^2-1)^2\sin^2(2\theta_1)\Big),     
\end{aligned} 
\]
\[
\begin{aligned}
{\beta(\theta_1)}=&b\Big(a^2-2a^4+\big(7-21a^2+8a^4\big)b^2+\big(7+8a^2\big)b^4-8b^6+\cos(2\theta_1)(a^2-b^2)\\&
(-1+2a^2-9b^2+8b^4)-8(a^2-1)^2 b^2\cos^2(2\theta_1)+2(a^2-b^2)(2a^2-b^2-1) \\&
\cos(4\theta_1)\sin^2(\theta_1)+8ab(b^2-1)^2\sin^2(2\theta_1)\Big),      
\end{aligned} 
\]
\[
\begin{aligned}
{\gamma(\theta_1)}=& -8b\Big(a^2(1-5b^2)-2b^4(b^2-2)+a^4(3b^2-1)-2(b^2-1)(b^2-a^2)\\&
 (1-a^2+b^2)\cos(2\theta_1)-(a^2-1)\big(a^2+(a^2-2)b^2\big)\cos^2(2\theta_1)+\\&
ab(b^2-1)^2\sin^2(2\theta_1)\Big).     
\end{aligned} 
\]
Also, for the function $\mu$ and constant $\xi$, we have
\begin{equation}\label{equation 4.40}
{\mu(\theta_1)}=(b^2-\sin^4(\theta_1)),\quad \xi=32(a^2-1)(b^2-1)(b^2-a^2). 
\end{equation}
These show that if $\theta_1=\theta_2$ then $\alpha(\theta_1)=\beta(\theta_1)$. Finding the maximum and minimum points of $Q$ is straight forward.  Set
\begin{equation}\label{equation 4.41} 
Q_1(\theta_1)=\frac{ b^2(\alpha(\theta_1)+\beta(\theta_1))-2b\;\sqrt{\alpha(\theta_1)}\sqrt{\beta(\theta_1)}\sin^2(\theta_1)+\gamma(\theta_1) \mu(\theta_1) }{\mu^2(\theta_1) \xi}    
\end{equation}
where $\theta_1\in [0,\pi]$. Then \begin{footnotesize} 
\begin{equation}\label{equation 4.42}
\begin{aligned}
\frac{\rm{d} Q_1}{\rm{d} \theta_1}(\theta_1)=& 2b^2\Big[\frac{\big(2+6b^2+4a^2(b+1)-4a(b+1)^2\big)\sin(2\theta_1)}{(a-b)(b-1)(a-1)\big(1+2b-\cos(2\theta_1)\big)^3}+\\&
\frac{\big(-1-2a^2(b+1)+2a(b+1)^2+b(-2+b-2b^2)\big)\sin(4\theta_1)}{(a-b)(b-1)(a-1)\big(1+2b-\cos(2\theta_1)\big)^3}\Big]=0\\&     
\end{aligned} 
\end{equation}
\end{footnotesize}
 There are five critical points on $[0,\pi]$, 
$$ \theta_1\in \{ 0,\frac{\pi}{2},\pi,\arctan\frac{y_1}{x_1},\pi-\arctan\frac{y_1}{x_1}, $$
where 
$$x_1=\frac{(b-1)\sqrt{b}}{\sqrt{1-2a+2a^2+2b-4ab+2a^2b-b^2-2ab^2+2b^3}},$$
and
$$y_1=\frac{\sqrt{b+1}\sqrt{1-2a+2a^2-2ab+b^2}}{\sqrt{1-2a+2a^2+2b-4ab+2a^2b-b^2-2ab^2+2b^3}}.$$
The values of function $Q_1$ at the critical points are 
$$Q_1(0)=0,\qquad\qquad Q_1(\frac{\pi}{2})=\frac{b}{b^2-1}>0,\qquad\qquad Q_1(\pi)=0,$$
and as $Q_1(\arctan\frac{y_1}{x_1})= Q_1(\pi-\arctan\frac{y_1}{x_1})$ and
$$ Q_1(\arctan\frac{y_1}{x_1})=  \frac{b(b-1)^3}{4(a-1)(b-a)(b+1)\big(1+a^2+b(b-1)-a(b+1)\big)},$$
these values are all non-negative. 
\subsubsection{Case 2: ($\theta_2=\pi-\theta_1$)} This is the case we worked out explicitly earlier at \S \ref{goodcase}. 

\subsection{The best rank-one direction.}
We have found the best rank-one direction, that which maximises the negative of the second derivative. 
We know that in the best direction, $\theta_2=\pi-\theta_1$ and $\alpha(\theta_1)=\beta(\theta_1)$. We have $\delta=r^2$ and $\eta=s^2$. Then 
$$\delta_1=\frac{b\Big(b-\sqrt{\frac{\beta}{\alpha}}\sqrt{b^2-\mu}\Big)}{\mu},  \qquad \eta_1=\frac{b\Big(b-\sqrt{\frac{\alpha}{\beta}}\sqrt{b^2-\mu}\Big)}{\mu}.$$
Since $\alpha(\theta_1)=\beta(\theta_1)$,  $r^2=s^2$, and $r=\pm s$. Without loss of generality we may assume that $r=s$, so, 
$$r=s=\frac{b\Big(b-\sqrt{\frac{\beta}{\alpha}}\sqrt{b^2-\mu}\Big)}{\mu}=\frac{\big(b^2-b\sin^2(\theta_1)\big)}{b^2-\sin^4(\theta_1)}=\frac{b\big(b-\sin^2(\theta_1)\big)}{(b-\sin^2(\theta_1))(b+\sin^2(\theta_1))}.$$
Hence,   $r=s=\frac{b}{(b+\sin^2(\theta_1))}$ and
\begin{equation}\label{equation 4.52}
\mathbf{u}^t=\Big(\frac{\sqrt{2b\sin(\theta_1)+\sin^2(\theta_1)}}{(b+\sin^2(\theta_1))},\;\frac{b\cos(\theta_1)}{(b+\sin^2(\theta_1))},\;\frac{b\sin(\theta_1)}{(b+\sin^2(\theta_1))}\Big),
\end{equation}
\begin{equation}\label{equation 4.53}
\mathbf{v}^t=\Big(\frac{\sqrt{2b\sin(\theta_1)+\sin^2(\theta_1)}}{(b+\sin^2(\theta_1))},\; -\frac{b\cos(\theta_1)}{(b+\sin^2(\theta_1))},\;\frac{b\sin(\theta_1)}{(b+\sin^2(\theta_1))}\Big),
\end{equation}
and with our previously computed values  $\mathbf{u}$ and $\mathbf{v}$ can be written with respect to $a$ and $b$ as $ \mathbf{u}=\big(u_1,u_2,u_3 \big),  \mathbf{v}=\big(u_1,-u_2,u_3 \big), $
where
\begin{equation}\label{uandv}
\begin{aligned}
&{u_1=\frac{(b-1)}{\sqrt{2(b+1)(1+a+a^2+b(a-1)+b^2)}}},\\[0.6 cm]
&{u_2=\sqrt{\frac{1+2a^2+b^2+2a(1+b)}{2(1+a+a^2+b(a-1)+b^2}}},\\[0.6 cm]
&{u_3=\frac{(b-1)\sqrt{b}}{\sqrt{2(b+1)(1+a+a^2+b(a-1)+b^2)}}}.
\end{aligned}
\end{equation}
The minimum of the function $Q$ is 
\begin{equation}\label{equation 4.55}
 Min({\rm diag}(1,a,b))=\frac{-b(b-1)^3}{4(a+1)(b+1)(a+b)(1+a+a^2+ab-b+b^2)}.
\end{equation}
\subsection{An  example: $A={\rm diag}(1,2,4)$}  Then $B_0=\mathbf{u} \otimes \mathbf{v} $,   
$ \mathbf{u}=\Big(\frac{1}{\sqrt{30}},\sqrt{\frac{5}{6}},\sqrt{\frac{2}{15}}\,\Big),\; 
\mathbf{v}=\Big(\frac{1}{\sqrt{30}},-\sqrt{\frac{5}{6}},\sqrt{\frac{2}{15}}\,\Big). $
$$A+tB_0=\begin{bmatrix} 
1+\frac{t}{30} & -\frac{t}{6} & \frac{t}{15}\\
 \frac{t}{6} & 2-\frac{5t}{6} & \frac{t}{3} \\
 \frac{t}{15} & -\frac{t}{3} & 4+\frac{2t}{15} \end{bmatrix}.$$
If $K=(A+tB_0)^t(A+tB_0)$, then the characteristic equation is 
$$-\lambda^3+\Big(21-\frac{11t}{5}+t^2\Big)\lambda^2+\Big(-84+50t-\frac{1121t^2}{100}\Big)\lambda+\Big(64-\frac{224t}{5}+\frac{196t^2}{25}\Big)=0.$$
Using Taylor's series the first three terms of the eigenvalues (on an interval around $0$) of the matrix $K$ are equal to
$$\begin{aligned}
&\lambda_1({t})=1+\frac{1}{15}{t}+\frac{1}{60}{t^2}+O({t^3}),\\&
\lambda_2({t})=4-\frac{10}{3}{t}+\frac{29}{36}{t^2}+O({t^3}),\\&
\lambda_3({t})=16+\frac{16}{15}{t}+\frac{8}{45}{t^2}+O({t^3}),
\end{aligned}$$
The linear distortion is 
$H(A+tB_0)=\sqrt{\frac{\lambda_3(t)}{\lambda_1(t)}}$.  
\begin{figure}[ht]
         \includegraphics[scale=0.5]{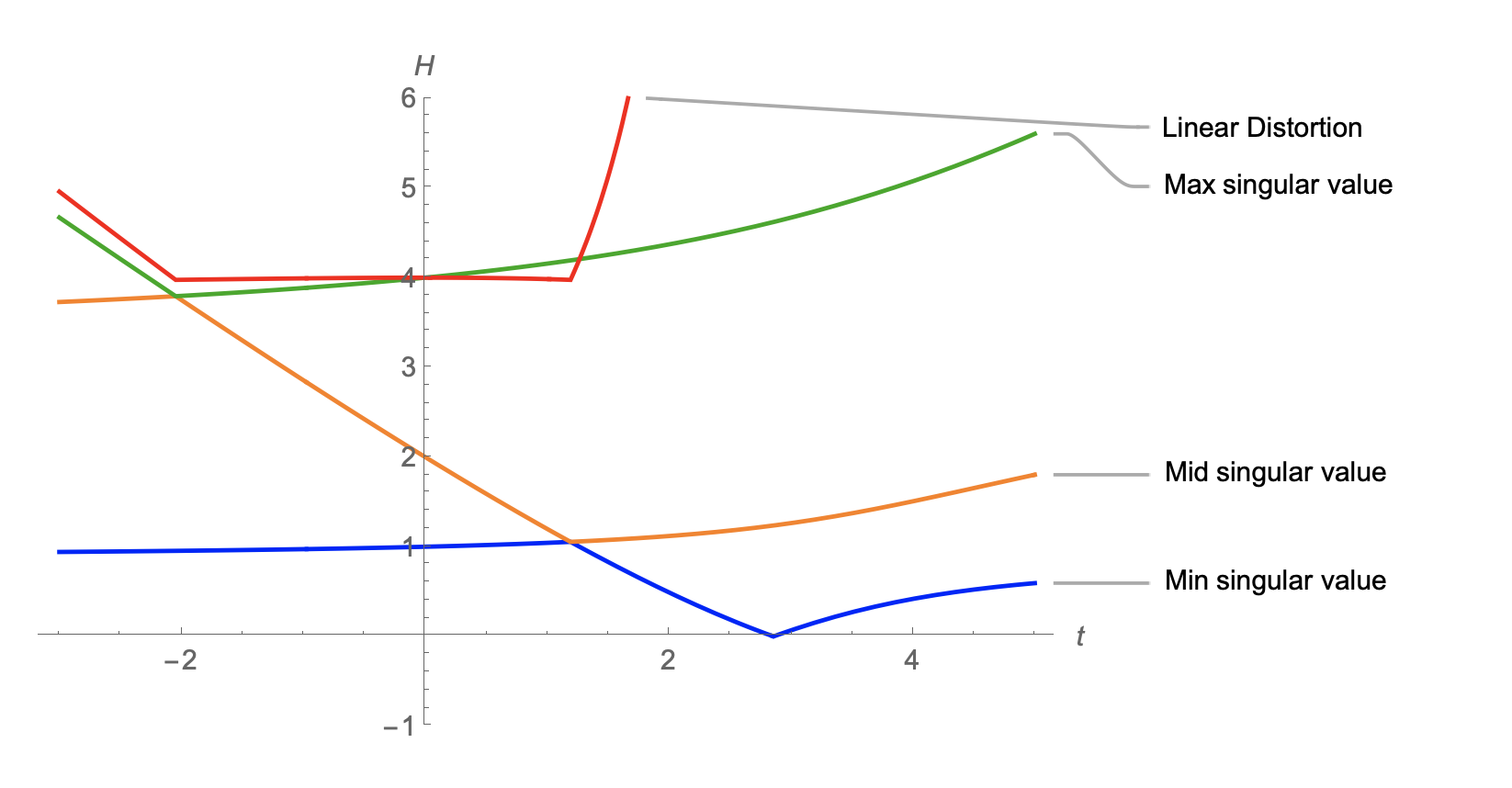}
      {{\bf Figure 6.} The eigenvalues of the matrix $K=(A+tB_0)^t(A+tB_0)$, when $A=(1,2,4)$  and $B_0$ is the optimal direction.}
         \label{fig: 4.24}     
\end{figure} 
$$H(A+tB_0)=\sqrt{\frac{\lambda_3(t)}{\lambda_1(t)}}=\sqrt{\frac{16+\frac{16t}{15}+\frac{8t^2}{45}}{1+\frac{t}{15}+\frac{t^2}{60}}}=4-\frac{1}{90}t^2+O(t^3). $$ 
$$\frac{\rm d}{\rm{d}t}\Big|_{t=0} H(A+t\, B)=  0, \;\;\; {\rm and} \;\;\; \frac{\rm d^2}{\rm dt^2}\Big|_{t=0} H(A+t\, B)=-\frac{1}{90}. $$
             Figure 6 also illustrates the loss of smoothness and concavity precisely where the eigenvalues cross.
             
\section{Solving Problem 2.}      Let  $A={\rm diag}(1,a,b)$ with $1<a<b$ and suppose $\mathbf{u}$ and  $\mathbf{v}$ in $ \mathbb{R}^{3} $ have $\|\mathbf{u}\|=\|\mathbf{v}\| =1$ so that $B_0=\mathbf{u} \otimes \mathbf{v}^t $   is a solution to Problem 1.  We want to find;
 \begin{enumerate}
\item  $\mathbf{t_+}>0$ and $\mathbf{t_-}<0$ of largest magnitude so that $ H(A+t\, B_0) $ is a smooth function of $t$ in the interval $\mathbf{t_-}<t<\mathbf{t_+}$,  and that
\item for all rank-one matrices $B$ with 
$$ \frac{\rm d}{\rm{d}t}\Big|_{t=0} H(A+t\, B)=  0, \;\;\; {\rm and} \;\;\; \frac{\rm d^2}{\rm dt^2}\Big|_{t=0} H(A+t\, B) < 0, $$ 
we have for all $t>0$
$$ H(A+t B) \geq \max \{ H(A+\mathbf{t_-} B_0),  H(A+\mathbf{t_+} B_0)\}. $$ 
\end{enumerate}
The values $\mathbf{t_-}$ and $\mathbf{t_+}$ we find are where the singular values of $A+tB_0$ cross. Thus we identify the discriminant of the eigenvalue equation for $H(A+tB_0)$ as it is the transverse crossing of the eigenvalues which implies $H$ loses smoothness. The two vectors $\mathbf{u}$ and $\mathbf{v}$ are given at (\ref{uandv}).  The characteristic equation $\det ((A+tB_0)^t(A+tB_0)-\lambda I)=0$ is:

\begin{footnotesize}
\[
\begin{aligned} 
&\det(K-\lambda I) =\\ & \frac{b^2 \Big(2 a^3 (b+1)+2 a^2 (b+1) (b-t+1)+2 a (b^3-4 bt+1)-(b^3+b^2+b+1) t\Big)^2}{4 (b+1)^2 \big(a^2+ab+a+b^2-b+1\big)^2}-\\&
\frac{1}{4 (b+1)^2 \big(a^2+ab+a+b^2-b+1\big)^2}\Big[4a^6(b+1)^2(1+b^2)+8a^5(b+1)^2(1+b^2)(1+b-t)+\\&
8b^6(t-2)t+t^2+4bt^2+4b^7t(t+2)-4b^3t(3t+4)+2b^4t(15t+4)+4b^2(1+2t+2t^2)+\\&
b^8(t^2+4)-4b^5(3t^2-2t-2)+4a^4(b+1)\Big(3+3b^5-3b^4(t-1)-3t+t^2+b(t^2-10t+3)+\\&
b^3(t^2-10t+7)+b^2(t^2-6t+7)\Big)-2a(b+1)\Big(2b^7t-(t-2)t-3b^4(t-2)-8bt^2+b^3(-8t^2+\\&
6t-4)-b^6(t^2+4)-4b^5(2t^2+1)-b^2(3t^2+4)\Big)+4a^3(b+1)\Big(2+2b^6+b^5(2-4t)-4t+t^2+\\&
b^4(t^2-8t+4)+2b^2(t^2-6t+2)+2b^3(3t^2-6t+4)+b(6t^2-8t+2)\Big)+a^2(b+1)\Big(4+4b^7-\\&
\end{aligned}
\]
\[
\begin{aligned} 
&8t+5t^2-4b^6(2t+1)+b^5(5t^2-28t+20)+2b^3(11t^2-36t+6)+b(21t^2-28t-4)+b^4(12+\\&
8t+21t^2)+b^2(20+8t+22t^2)\Big)\Big]\,{\lambda}+\frac{1}{(b+1)(1+a+a^2-b+ab+b^2)}\Big[1+b^5+a^4(b+1)+\\&
a^3(b+1)(1+b-2t)+b^3(t-1)^2+ t-2bt+b^4t+t^2+b^2(1+2b)+a(1+b)(1+b+b^2+\\&
b^3-t-b^2t+t^2+bt^2)+a^2(b+1)\big(2+2b^2-2t+t^2-b(1+2t)\big) \Big]\,{\lambda^2} -{\lambda^3}=0.
\end{aligned}
\]
\end{footnotesize}
While degree three in $\lambda$,  this polynomial is only quadratic in $t$ -- a consequence of the Jacobi identity for determinants.
This characteristic equation must have three nonnegative real roots. To avoid solving these equations we use an algebraic trick: A root in $t$ to the above equations implies a repeated root of the characteristic equation.

\begin{theorem} The characteristic polynomial
\[ \det ((A+tB_0)^t(A+tB_0)-\lambda I)=0 \]
has exactly two real roots as a polynomial in $t$.
\end{theorem}
We write out this equation in $t$. Remarkably the discriminant equation of degree $8$ in $t$ has a simple repeated quadratic factor $P(t)$ equal to 
 \[
\begin{aligned}
&\big(-1-a-2a^2-4b-7ab-4a^2b+2b^2-7ab^2-2a^2b^2-4b^3-ab^3-b^4\big)\;{t^2}+\\&
\big(-1+2a+a^2+4a^3-6b-4ab+7a^2b+8a^3b-b^2-12ab^2+7a^2b^2+4a^3b^2-\\&
b^3-4ab^3+a^2b^3-6b^4+2ab^4-b^5\big)\;{t}-2\big(-a+a^4+b-ab-2a^2b\\&+2a^4b+b^2+
2ab^2-4a^2b^2+a^4b^2+2ab^3-2a^2b^3+b^4-ab^4+b^5-ab^5\big)
\end{aligned}
\]
The quartic remainder has discriminant in $t$ equal to
\begin{footnotesize}
\[ \begin{aligned}
&256 (a+1)^2 (b-1)^{12} b (b+1)^7 (a+b)^2 \left(a^2+a b+a+b^2-b+1\right)^6 \left(2 a^2+2 a (b+1)+b^2+1\right) \\ &\cdot \left(27 a^6 \left(b^3+b^2+b+1\right)+54 a^5 (b+1)^2 \left(b^2+b+1\right)+9 a^4 \left(7 b^5+27 b^4+56 b^3+56 b^2+27 b+7\right)\right.\\&\left. +3 a^2 +4 a^3 \left(b^2+b+1\right)^2 \left(11 b^2+38 b+11\right)+\left(7 b^7+46 b^6+99 b^5+118 b^4+118 b^3+99 b^2+46 b+7\right)\right.\\&\left.+6 a (b+1)^2 \left(b^6+6 b^5+3 b^4+7 b^3+3 b^2+6 b+1\right)\right.\\&\left.+b^9+6 b^8+18 b^7+8 b^6+21 b^5+21 b^4+8 b^3+18 b^2+6 b+1\right)^3 \end{aligned}
\]
\end{footnotesize}
This is strictly positive.  The discriminant of the second derivative is negative, so has no real roots,  and it follows that this quartic does not have four real roots and so it has none.  

\begin{corollary} The regular branches of the eigenvalues of $(A+tB_0)^t(A+tB_0)$ cross twice.
\end{corollary} 
In fact these crossings are transverse and so $H(A+tB_0)$ will lose smoothness there.  We do not need this result, but it can be proved by a lengthy calculation from what follows as we calculate these crossing points.
We have
 
\[
\rm{Disc.}[\det ((A+tB_0)^t(A+tB_0)-\lambda I)]=\frac{P(t) ^2 R(t)}{16(b+1)^5(1+a+a^2-b+ab+b^2)^5}\; ,
\]
with $R(t)>0$. The discriminant of  $P(t)$ is 
\[ 
\begin{aligned}
&\Delta= (b^2-1)^2\Big[1+4b+6b^3+4b^5+b^6+4a(b+1)\big(b^4+3b^3+3b+1\big)\\&+4a^3(b+1)
\big(3+(2+3b)\big)+2a^2(1+b+b^2)\big(5+b(6+5b)\big)+a^4\big(9+b(9b-2)\big)\Big],
\end{aligned}
\]
which is obviously positive. Thus $P(t)$ has two real roots and examining the coefficients shows they have different signs.   Let us denote these roots as $\mathbf{t_+}>0$ and $\mathbf{t_-}<0$. We calculate that
\[ 
\begin{aligned}
&\mathbf{t_+}=\frac{G_1(a,b)-(b^2-1)\sqrt{J_1(a,b)}}{2(-1-a-2a^2-4b-7ab-4a^2b+2b^2-7ab^2-2a^2b^2-4b^3-ab^3-b^4)},\\&
\mathbf{t_-}=\frac{G_1(a,b)+(b^2-1)\sqrt{J_1(a,b)}}{2(-1-a-2a^2-4b-7ab-4a^2b+2b^2-7ab^2-2a^2b^2-4b^3-ab^3-b^4)},
\end{aligned}
\]

where, 
\[
\begin{aligned}
G_1(a,b)=&1-2a-a^2-4a^3+6b+4ab-7a^2b-8a^3b+b^2+12ab^2-7a^2b^2-\\&
4a^3b^2+b^3+4ab^3-a^2b^3+6b^4-2ab^4+b^5,
\end{aligned}
\]
and
\[
\begin{aligned}
J_1(a,b)=&1+4a+10a^2+12a^3+9a^4+4b+16ab+22a^2b+20a^3b-2a^4b+\\&
12ab^2+32a^2b^2+20a^3b^2+9a^4b^2+6b^3+12ab^3+22a^2b^3+12a^3b^3+\\&
16ab^4+10a^2b^4+4b^5+4ab^5+b^6.
\end{aligned}
\]
Write the characteristic equation as
\begin{equation}\label{equation 4.63}
D_1{\lambda^3}+C_1{\lambda^2}+B_1{\lambda}+A_1=0,
\end{equation}
where
\begin{footnotesize}
$$\begin{aligned}
&D_1=-4(b+1)^2\big(1+a+a^2+(a-1)b+b^2 \big)^2,\\&
C_1=4(b+1)\big(1+a+a^2+(a-1)b+b^2 \big)\Big(1+b^2+b^3+b^5+a^4(b+1)+a^3(b+1)\\&
\qquad\:\:\,(1+b-2t)+(b-1)^2(b^2+1)t+(b^3+1)t^2+a^2(b+1)\big(2-b+2b^2-2(b+1)t+\\&
\qquad\:\:\,t^2\big)+a(b+1)\big(1+b+b^2+b^3-(b^2+1)t+(b+1)t^2\big)\Big),\\[0.2cm]
\end{aligned}$$
$$\begin{aligned}
&B_1=-4a^6(b+1)^2(b^2+1)-4(b^4+b)^2-8a^5(1+b)^2(b^2+1)(1+b-t)-8b^2(b-1)^2\\&
\qquad\:\:\, (b^3+1)t-\Big[1+b\Big(4+b\big[8+b\big(-12+b(30+b(-12+b(8+4(4+b)))) \big)\big]\Big)\Big]+\\&
\qquad\:\:\, 4a^4(b+1)\Big(-(b+1)(3+7b^2+3b^4)+(b+3)(3b+1)(b^2+1)t-(b+1)(b^2+1)t^2\Big)+  \\&
\qquad\:\:\, 4a^3(b+1)\Big(-2\big(1+b+b^2(b+1)(2+2b+b^3)\big)+4(b+1)(b^2+1)(1+b+b^2)t - \\&
\qquad\:\:\,(b^2+1)(1+b(6+b))t^2\Big)+a^2(b+1)\Big[-4(b+1)(b^2+1)\Big(1+b\big(-2+b(6+    \\&
\qquad\:\:\, (-2+b)b)\big)\Big)+4\Big(2+b\Big(7+b\Big[-2+b\Big(18+b\big(-2+b(7+2b)  \big)\Big)\Big]\Big)\Big)t-(1+b)\\&
\qquad\:\:\, \Big(5+b\Big(16+b\big(9+b(16+5b)  \big)\Big)\Big)t^2\Big]+ 2a(b+1)\Big[-4(b^2+b^3+b^5+b^6)+2t\\&
\qquad\:\:\, \big(1+b^3(3+3b+b^4)\big)-\Big(1+b\Big(8+b\Big(3+b\big(8+b(3+b(8+b))   \big)\Big)\Big)\Big)t^2\Big] ,\\&
A_1=b^2\Big(2a^3(b+1)+2a^2(b+1)(1+b-t)-(b+1)(b^2+1)t+2a(1+b^3-4bt) \Big)^2.
\end{aligned}$$
\end{footnotesize}
The roots of (\ref{equation 4.63}) are 
\begin{equation}\label{equation 4.64}
\begin{aligned}
&{\lambda_1}=Z(a,b,t)-\frac{2^{\frac{1}{3}}\times Y(a,b,t)}{3\times X(a,b,t)}+\frac{X(a,b,t)}{3\times2^{\frac{1}{3}}\times D_1}  ,\\&
{\lambda_2}= Z(a,b,t)+\frac{\frac{1-i\sqrt{3}}{2}\times 2^{\frac{1}{3}}\times Y(a,b,t)}{3\times X(a,b,t)}-\frac{\frac{1+i\sqrt{3}}{2}\times X(a,b,t)}{3\times2^{\frac{1}{3}}\times D_1}   ,\\&
{\lambda_3}=Z(a,b,t)+\frac{\frac{1+i\sqrt{3}}{2}\times 2^{\frac{1}{3}}\times Y(a,b,t)}{3\times X(a,b,t)}-\frac{\frac{1-i\sqrt{3}}{2}\times X(a,b,t)}{3\times2^{\frac{1}{3}}\times D_1}  
\end{aligned}
\end{equation}
where 
\begin{eqnarray*} X(a,b,t)&=& \Big( -2C^3_1+9B_1C_1D_1-27A_1D_1 \\ && +\sqrt{-4\big(C^2_1-3B_1D_1 \big)^3+\big(2C^3_1-9B_1C_1D_1+27A_1D^2_1  \big)^2}\;\Big)^{\frac{1}{3}},\\
Y(a,b,t)&=&\frac{(-C^2_1+3B_1D_1)}{D_1}\\ Z(a,b,t)&=&-\frac{C_1}{3D_1}.
\end{eqnarray*} 
The linear distortion of the matrix $A+tB_0$ is:

\begin{equation}\label{equation 4.65}
H(A+tB_0)=\sqrt{\frac{\lambda_3}{\lambda_1}}.
\end{equation}
 
We are now in a position to solve our problem.  Given $1<a<b$ we find $B_0$, as at (\ref{uandv}) and calculate  $H(A+\mathbf{t_+}B_0)$ and $H(A+\mathbf{t_-}B_0)$. Both of these are less than $b=H(A)$. 
 

Following on from our earlier example for $A=\mathrm{{\rm diag}}(1,2,4)$, we compute that 
$$H(A)=4, \;\; \mathbf{t_+}=1.19219, \;\;\mathbf{t_-}=-2.04584,$$ $$H(A+\mathbf{t_+}B_0)=3.97539, \;\; H(A+\mathbf{t_-}B_0)=3.97539.$$
So, for all $t\in(-2.04484,\;1.19219)$, we have $H(A+tB_0)<H(A)$.  This example is remarkable in that $H(A+\mathbf{t_+}B_0)= H(A+\mathbf{t_-}B_0)$.  In general this is not true,  but appear to be the case in the situation we believe is extremal.

\section{Iwanec's construction} Here we briefly sketch this construction.
 We have $ B=\mathbf{u} \otimes \mathbf{v} $. Consider $ \mathbf{T_{0}}(x)= A x $, and for $ \nu=1,2,\ldots $, we define a sequence $ \{ \mathbf{T_{\nu}} \}_{\nu=1}^{\infty} $ by equation
$$ \mathbf{T_{\nu}}(x)=Ax \frac{1}{\nu}\,h\,(\nu \mathbf{u}\cdot x)\mathbf{v} $$
where $h$ is a periodic piecewise linear function on the real line illustrated in Figure 7.
\begin{center} \includegraphics[scale=0.4]{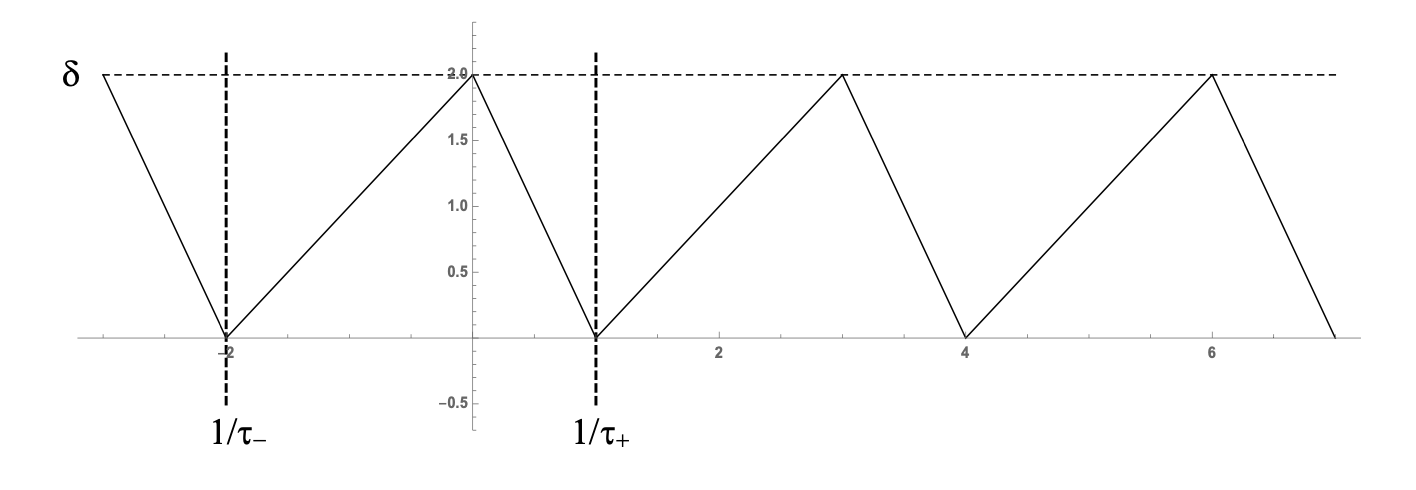}
\end{center} 
 {\bf  Figure 7.} {The saw-tooth function}   
Given $\mathbf{t_-}<0$ and $\mathbf{t_+}>0$ of the lemma we define
$$ h(r)=
\begin{cases}
\mathbf{t_-}r               &\qquad \frac{i}{\mathbf{t_+}}-\frac{i-1}{\mathbf{t_-}} \leq r \leq \frac{i}{\mathbf{t_+}}-\frac{i}{\mathbf{t_-}}, \\
\mathbf{t_+}r              &\qquad \frac{i}{\mathbf{t_+}}-\frac{i}{\mathbf{t_-}} \leq r \leq \frac{i+1}{\mathbf{t_+}}-\frac{i}{\mathbf{t_-}},
\end{cases}
$$
for any integer $i$. Then we can extend $h$ to the entire line (saw-tooth function). The function $h$ is a bounded Lipschitz function whose derivative assumes only the two values $\mathbf{t_-}$ and $\mathbf{t_+}$. The sequence $ \{ \mathbf{T_{\nu}} \}_{\nu=1}^{\infty} $ converges uniformly to $ \mathbf{T_{0}} $. The derivative of $ \mathbf{T_{\nu}} $ also assumes only two values, which are independent of $\nu$, apart from countable set of points where it is not defined. 
$$ D\mathbf{T_{\nu}} = A+h'(\nu \mathbf{u}\cdot x)\mathbf{u} \otimes \mathbf{v}=A+h'(\nu \langle \mathbf{u},x\rangle)B \in \{ A-\mathbf{t_-}B,A+\mathbf{t_+}B\}. $$
In either case, the linear distortion of $ D\mathbf{T_{\nu}}(x) $ is equal to $H$,
\[ H =\max\{H(A+\mathbf{t_-}B_0),H(A+\mathbf{t_+}B_0) \} \]  
Iwaniec's argument now proves our first theorem. We remark that the sequence $\{h'(\nu \langle\mathbf{u}, x\rangle)\}_{\nu=1}^{\infty}$ converges weakly in $L^{\infty}(\mathbb{R}^3)$ to $0$ as $\nu \longrightarrow \infty$, but not pointwise almost everywhere.

\subsection{Higher dimensions.}

It is very hard to identify an optimal rank-one direction in $n$-dimensions, because the characteristic equations of $A+tB$ is so complicated. However it is clear that the above arguments work in higher dimensions when we extend the matrices $A$ and $B$ by the rules (see \cite{I})
$$A=\begin{bmatrix} {1} & 0 & 0 & 0 & \cdots & 0 \\ 0 & {a} & 0 & 0 &  \cdots & 0 \\ 0 & 0 & {b} & 0 &\cdots & 0\\ 0 & 0 & 0 & {a} & \cdots & 0 \\  \vdots & \vdots & \vdots & \vdots &   \ddots & \vdots \\0 & 0 & 0 & 0 &  \cdots & {a} \end{bmatrix}\quad {\rm
and} \quad  B=\begin{bmatrix}  B_0 & 0 & \cdots & 0\\ 0   & 0 &  \cdots & 0   \\ \vdots & \vdots & \ddots & \vdots   \\0 & 0  & \cdots & 0 \end{bmatrix}. $$

 \section{The bounds.}

For $A={\rm diag}(1,a,b)$ we have everything written explicitly in terms of $a$ and $b$,  however complicated, and so we may explore where the largest jump and the number $\sqrt{2}$ comes from.  We consider $A={\rm diag}(1,c,f(c))$,  identify the best direction and maximal jump as functions of $c$.  These are illustrated for a variety of choices of $f$ below in Figures 8 and 9.  These calculations led directly to the conjecture that this jump is largest when $A={\rm diag}(1,c,c^2)$ and $c$ is large,  remarkably the example that Iwaniec considered though he did not have the optimal direction.
 \begin{multicols}{2}
      \includegraphics[scale=0.42]{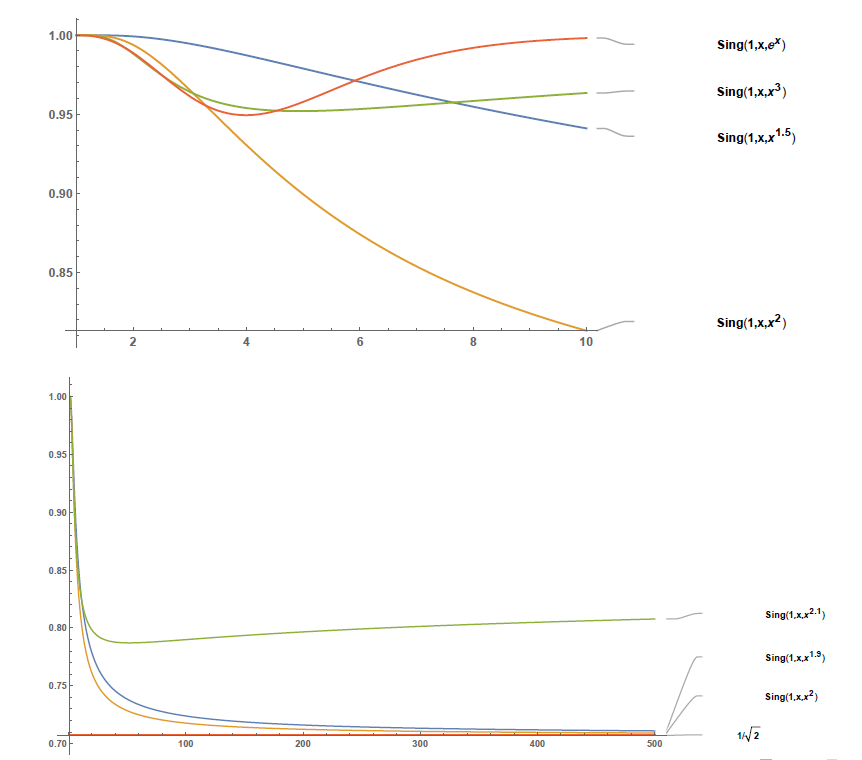} \par  
\quad \includegraphics[scale=0.33]{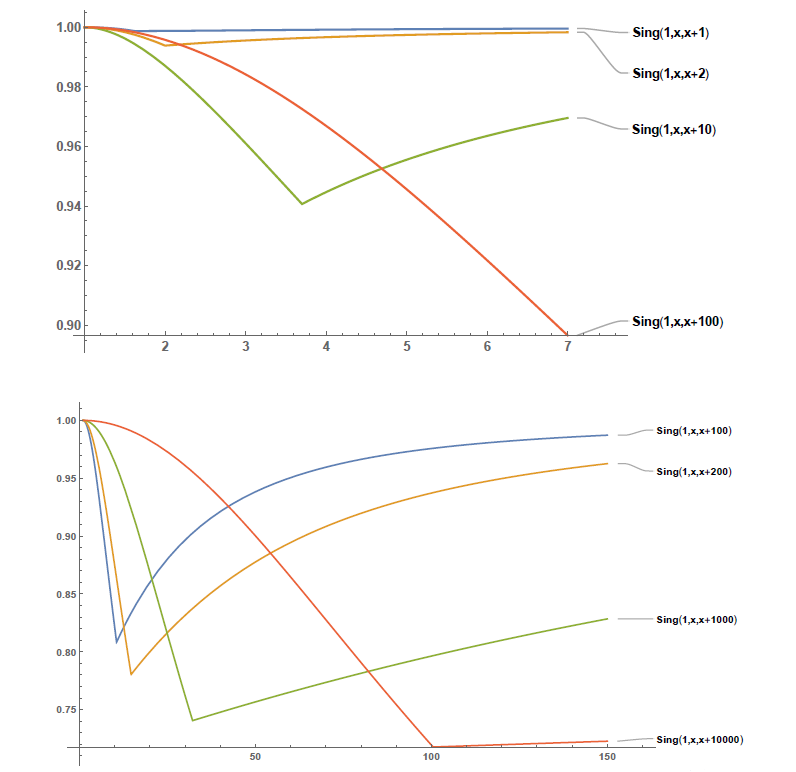}\par  
\end{multicols}
\noindent {\bf Figure 8} The jump: $1<<a<<b$ $\quad\quad $ {\bf Figure 9} The jump $|a-b|$ bounded. 

The optimal direction is
\begin{small}\[  B_0=\left(
\begin{array}{ccc}
 \frac{1}{c^2+1}-\frac{1}{2 (c-1) c+2} & -\frac{\left(c^2-1\right) \sqrt{\frac{c}{(c-1) c+1}+1}}{2 \sqrt{(c+1)^2 ((c-1) c+1) \left(c^2+1\right)}} & \frac{(c-1)^2 \sqrt{c^2}}{2 ((c-1) c+1) \left(c^2+1\right)} \\
 \frac{\left(c^2-1\right) \sqrt{\frac{c}{(c-1) c+1}+1}}{2 \sqrt{(c+1)^2 ((c-1) c+1) \left(c^2+1\right)}} & -\frac{c^2+1}{2 (c-1) c+2} & \frac{\sqrt{c^2} \left(c^2-1\right) \sqrt{\frac{c}{(c-1) c+1}+1}}{2 \sqrt{\left(c^2+1\right) \left(c^4+c^3+c+1\right)}} \\
 \frac{(c-1)^2 \sqrt{c^2}}{2 ((c-1) c+1) \left(c^2+1\right)} & -\frac{\sqrt{c^2} \left(c^2-1\right) \sqrt{\frac{c}{(c-1) c+1}+1}}{2 \sqrt{(c+1)^2 ((c-1) c+1) \left(c^2+1\right)}} & \frac{(c-1)^2 c^2}{2 ((c-1) c+1) \left(c^2+1\right)} \\
\end{array}
\right) \]
\end{small}
\begin{scriptsize}
\begin{eqnarray*}
\mathbf{t_+}& = & \frac{(c-1) \left(c^2+1\right) \left(-c^6+2 c^5-5 c^4+5 c^2+(c+1)^2 \sqrt{\left(c^2+1\right) \left(\left(c \left(c^3+7 c-8\right)+7\right) c^2+1\right)}-2 c+1\right)}{2 \left(\left(c^5+6 c^3+c^2+c+6\right) c^2+1\right)}\\
\mathbf{t_-}& = & -\frac{(c-1) \left(c^2+1\right) \left(c^6-2 c^5+5 c^4-5 c^2+(c+1)^2 \sqrt{\left(c^2+1\right) \left(\left(c \left(c^3+7 c-8\right)+7\right) c^2+1\right)}+2 c-1\right)}{2 \left(\left(c^5+6 c^3+c^2+c+6\right) c^2+1\right)}
\end{eqnarray*}
\end{scriptsize}
With these values one can compute the singular values of $A+\mathbf{t_\pm}B_0$ explicitly,  and then $H(A+\mathbf{t_\pm}B_0)$.  The formula for $H(A+\mathbf{t_\pm}B_0)$ runs over a few pages and we do not reproduce it here.  However is involves only Laurent polynomials in $c$ of modest degree $18$ and the square roots of polynomials of degree $8$ and a cube root.  As such the limit as $c\to\infty$ can be directly computed using the usual bag of tricks,  we checked our results also with Mathematica which returns the limit in a few minutes.  With $A$, $B_0$ and $\mathbf{t_\pm}$ as above we found
\begin{equation}
\lim_{c\to \infty} \frac{H(A+\mathbf{t_\pm}B_0)}{c^2} =\frac{1}{\sqrt{2}}
\end{equation}
and this is our second theorem.

Institute for Advanced Study, 
Massey University,  Auckland,
New Zealand.
\newline
email: g.j.martin@massey.ac.nz 
\end{document}